\documentclass{ifacconf}

\usepackage{graphicx}      
\usepackage{natbib}        
\usepackage{times,latexsym,amssymb,amsmath,theorem,epsfig}

\newtheorem{theorem}{Theorem}[section]

\newtheorem{definition}[theorem]{Definition}
{\theorembodyfont{\rmfamily} \newtheorem{example}[theorem]{Example}}

\newtheorem{problem}[theorem]{Problem}

\begin{document}
\newcommand{\diag}{\mbox{${\rm diag}$}}
\newcommand{\dnn}{\mbox{$\mathbb{R}_{\rm diag}^{n \times n}$}}
\newcommand{\dposnn}{\mbox{$\mathbb{R}_{{\rm diag},+}^{n \times n}$}}
\newcommand{\dsposnn}{\mbox{$\mathbb{R}_{{\rm diag},s+}^{n \times n}$}}
\newcommand{\argmin}{\mbox{${\rm argmin}$}}
\newcommand{\Table}{\mbox{$\rm table$}}
\newcommand{\gmon}{\mbox{${\rm G_{mon}}$}}
\newcommand{\orderlexg}{\mbox{$>_{{\rm lex}}$}}
\newcommand{\zn}{\mbox{$\mathbb{Z}_n$}}
\newcommand{\card}{\mbox{$\rm card$}}
\newcommand{\matrixdiagtrunc}{\mbox{$\rm mdtrunc$}}
\newcommand{\svdtruncation}{\mbox{${\rm SVDtrunc}$}}
\begin{frontmatter}

\title{
An Algorithm for System Identification of 
a Discrete-Time Polynomial System without Inputs \\
-- Extended Version
\thanksref{footnoteinfo}} 

\thanks[footnoteinfo]{
Jana N\v{e}mcov\'{a} thanks 
the national agency of the Czech Republic 
for financial support by
GA\v{C}R Project 13-16764P.
}

\author[First]{Jana N\v{e}mcov\'{a},} 
\author[Second]{Mih\'{a}ly Petreczky,} 
\author[Third]{Jan H. van Schuppen} 

\address[First]{
Department of Mathematics, Institute of Chemical Technology, \\
Prague, Czech Republic
(e-mail: Jana.Nemcova@vscht.cz)}
\address[Second]{
Dept. of Comp. Science \& Automatic Control,
\'{E}coles des Mines de Douai,\\
Douai, France 
(e-mail: mihaly.petreczky@mines-douai.fr)}
\address[Third]{
Van Schuppen Control Research, 
Amsterdam, The Netherlands\\
(e-mail: jan.h.van.schuppen@xs4all.nl)}
\begin{abstract}                
A subalgebraic approximation algorithm 
is proposed to estimate from a set of time series
the parameters of the observer representation
of a discrete-time polynomial system without inputs 
which can generate an approximation
of the observed time series.
A major step of the algorithm is to construct
a set of generators 
for the polynomial function
from the past outputs to the future outputs.
For this singular value decompositions and
polynomial factorizations are used.
An example is provided.
\end{abstract}


\end{frontmatter}
%
\section{Introduction}
The system identification of polynomial systems is
motivated by the need for models of biochemical
reaction systems in the life sciences.
Also in the area of control engineering and
of economics there is a need to determine
parameter values of such control systems from data.
As far as the authors have been able to determine
there is no satisfactory algorithm for the general problem
of determining the parameter values of these systems.
\par
The problem of the paper is to determine
a system in the class of discrete-time 
polynomial systems without inputs 
in the form of an observer realization
such that it produces for each time series of outputs 
a predicted time series
which is a reasonable approximation
of the supplied output time series. 
\par
The relevant literature on polynomial systems and their
system identification is briefly summarized.
At the time this paper is written
there are available results on 
the realization theory of 
polynomial and of rational systems
see \cite{sontag:1979:phdthesis},
\cite{bartosiewicz:1988},
\cite{nemcova:schuppen:2009},\\
\cite{nemcova:schuppen:2010}.
The problem of structural identifiability
of polynomial and of rational systems was solved
by J. N\v{e}mcov\'{a} in \cite{nemcova:2010:mb}.
Earlier papers of the authors include
\cite{nemcova:schuppen:2009:sysid}, 
\cite{nemcova:petreczky:schuppen:2012:sysid}.
Various aspects of system identification of rational systems 
are also discussed in \\
\cite{boulier:lemaire:2009:sysid},
\cite{gevers:bazanella:coutinho:dasgoupta:2013}, \\
\cite{bazanella:gevers:coutinho:rui:2014}.
\par
The algorithm proposed in this paper 
determines a polynomial system
in the form of an observer,
thus driven by the available output.
It will be proven using system theory that
an observer polynomial system is equivalent with
the conditions: 
(1) the state of the observer at any time 
is a polynomial function of the past outputs; 
(2) the future outputs are a polynomial function only
of the current state; and
(3) the next state is a polynomial function only 
of the current state and the current output.
The main step of the algorithm is 
to construct an approximate generator set 
for the polynomial equation
from the past outputs to the future outputs.
The subalgebra generated by this generator set 
is then an approximate subalgebra of 
the algebra of the function from past outputs
to future outputs.
\section{Problem Formulation}
System identification of a polynomial and rational system
is motivated by the occurence of these systems in
engineering (satellite orientation problems),
biochemical reaction networks (mass action kinetics), and
economics (products of prices and quantities).
In this paper the authors focus on polynomial systems.
For the extension to rational systems and
to systems with inputs 
there is insufficient space in this short paper.
\par
System identification is a research area
that addresses the problem of
how to go from observational data
to a system with its parameter values.
The following procedure is often used:
(1) {\em Modeling.}
Model the phenomenon as a control system
as understood in system theory;
(2) {\em Data collection.}
Collect data in the form of a time series
of the phenomenon to be identified;
(3) {\em Structural identifiability.}
Determine whether 
the selected parametrization of the subclass of systems 
is structurally identifiable
and, if not, modify the system subclass such
that the system subclass is structurally identifiable;
(4) {\em Approximation.}
Determine an algorithm to compute 
the parameters of a system in the considered subclass
from the observation data which is a reasonable
approximation according to an approximation criterion;
(5) {\em Complexity estimation.}
Determine that subclass of systems which achieves
a reasonable value for the approximation criterion
and which minimizes the complexity. 
This paper addresses Step (4).
\par
For the approximation problem of system identification
there are two main approaches:
(1) Minimization of an approximation criterion
over the considered subclass of systems.
This approach, though often used, suffers
from the problem that the criterion 
is a nonconvex function of the parameters
which makes the minimization approach 
practically almost impossible.
(2) An algebraic method based on realization theory
of system theory.
The method is explained in the next section.
\par
For the remainder of this paper
the reader is expected to have read the notation
and terminology of the appendices. 
\begin{definition}\label{def:timeseries}
A {\em time series}
is a collection of real vectors 
denoted by 
its dimension $d_y \in \mathbb{Z}_+$,
its length in time steps $t_1 \in \mathbb{Z}_+$, and
$\{ \overline{y}(t) \in \mathbb{R}^{d_y}, ~ t \in \{ 1,2, \ldots, t_1\} \subset \mathbb{Z} \}$.
Often there is available a finite set of time series.
\end{definition}
\begin{definition}\label{def:observersystem}
The system class considered is that of
a discrete-time polynomial system without inputs
in observer representation,
which produces $y(t|t-1)$,
a prediction of the next output based on past outputs.
\begin{eqnarray*}
      x(t+1) 
  & = & f_o(x(t),y(t)), ~ x(0) = x_0, \\
      y(t|t-1)
  & = & h_o(x(t)), \\
      T
  & = & \{ 0,1,2, \ldots, t_1\} \subset \mathbb{N}, ~~
        X = \mathbb{R}^n, ~
        Y = \mathbb{R}^{d_y}, \\
  &   & f_o: X \times Y \rightarrow X, ~
        h_o: X \rightarrow Y, 
\end{eqnarray*}
where $f_o$ and $h_o$ are both polynomial functions.
\end{definition}
\begin{problem}\label{problem:siapproximation}
Consider a finite set of time series and
the subclass of control systems
of Def. \ref{def:observersystem}.
The problem is to determine 
an observer system in the subclass considered, 
specified by its parameter values,
such that the estimated observer system,
when supplied with the output of the time series,
produces a one-step prediction of the output
which prediction time series
is close to the observed time series.
\end{problem}
\section{The Approach}\label{sec:approach}
The subalgebraic approach to system idenfication
of polynomial systems is based on the
subspace identification algorithm 
of Gaussian systems. 
The approach was initially proposed by H. Akaike
based on contacts with R.E. Kalman,
has been developed for Gaussian systems,
and is known as the subspace identification algorithm,
see for references \cite{vanoverschee:demoor:1996}.
The subalgebraic ap\-pro\-ach to polynomial systems
is based on realization theory of nonlinear
systems in particular 
of bilinear, of polynomial and of rational systems.
See for references on realization theory
those mentioned in the previous section and
\cite{fliess:1990:forummath}, 
\cite{isidori:1973},
\cite{vidal:2008}.
No approximation criterion is used in the paper.
The principle of the subalgebraic method 
is explained with the following theorem.
\begin{theorem}\label{theorem:systemfactorization}
Consider a discrete-time polynomial system\\
with\-out inputs
in its observer realization 
with as output the one-step prediction,
\begin{eqnarray}
      x_o(t+1)
  & = & f_o(x_o(t),y(t)), ~ x_o(0) = x_{o,0}, 
          \label{eq:observerpolsys} \\
      y(t|t-1)
  & = & h_o(x_o(t)), ~
        f_o, ~ h_o, ~ 
	\mbox{polynomial maps.} \label{eq:observeroutputeq}
\end{eqnarray}
The observer system representation
(\ref{eq:observerpolsys},\ref{eq:observeroutputeq})
may be transformed to 
the following set of polynomial functions
assuming that the observer is a true observer,
hence the predictions equal $y(t) = y(t|t-1)$, and
in terms of $(y^+(t),~ y^-(t))$ 
as defined in equation (\ref{eq:yplusymin}),
\begin{eqnarray}
      x_o(t)
  & = & g_{io}(x_{o,0},y^-(t)), \label{eq:x0ypast} \\
      y^+(t)
  & = & h_{io}(x_o(t))
        = h_{io}(g_{io}(x_{o,0},y^-(t))),
        \label{eq:yfuturex0}  \\
      x_o(t+1)
  & = & f_o(x_o(t),y(t)), ~
        f_o, ~ g_{io}, ~ h_{io}, ~ \mbox{polyn.}  
        \label{eq:dynamics}
\end{eqnarray}
\end{theorem}
Note that 
by the equations (\ref{eq:x0ypast},\ref{eq:yfuturex0})
the current state is 
a polynomial of the components of the past outputs and
the future outputs are polynomial functions
of the components of the current state; and
by equation (\ref{eq:dynamics})
the next state $x_o(t+1)$ is a polynomial map
in the tuple $(x_o(t), y(t))$
of the current state and the current output.
If one considers the initial state as a constant
then the equations 
(\ref{eq:x0ypast},\ref{eq:yfuturex0})
define for any time a polynomial function 
from the past outputs to the future outputs.
\par
The algorithm for the subalgebraic approximation
is based on the above theorem
and consists of the steps:
\begin{enumerate}
\item
Compute a state vector $x(t)$ in terms of 
a polynomial function of past outputs
such that 
the future outputs are a polynomial function of it.
In terms of formulas,
\begin{eqnarray}
      y^+(t)
  & \approx & f_{io}(y^-(t)) 
        = h_{io}(g_{io}(y^-(t))) 
        = h_{io}(x(t)), 
        \label{eq:iofactorization} \\
      x(t)
  & = & g_{io}(y^-(t)), ~~ 
        f_{io}, ~ g_{io}, ~ h_{io} ~
        \mbox{polynomial functions.} \nonumber
\end{eqnarray}
\item
Compute the polynomial system dynamics 
\begin{eqnarray}
      x(t+1)
  & \approx & f_o(x(t),y(t)). 
        \label{eq:systemdynamicspolynomial}
\end{eqnarray}
\end{enumerate}
The main task of the algorithm
is to compute a set of generators 
of the polynomial map $f_{io}$
from past outputs to future outputs.
The complexity of the computations is limited 
by several steps of the algorithm.
\section{The Algorithm}\label{sec:algorithm}
\begin{definition}\label{alg:siapolsysdtnoinput}
{\em The subalgebraic approximation algorithm}
for system identification 
of discrete-time polynomial systems.
Data: 
A time series of outputs with the notations:
the dimension of the output $d_y \in \mathbb{Z}_+$,
the length of the time series $t_1 \in \mathbb{Z}_+$, 
the number of time series $s \in \mathbb{Z}_+$, and
finally the time series matrix
$Y_{ts} \in \mathbb{R}^{t_1 \times dy \times s}$. 
The parameters of the algorithm are
$r_1, ~ r_2, ~ r_3, ~ r_4 \in (0,1) \subset \mathbb{R}$ and
$(t_{\min}^+, ~ t_{\min}^-), ~ (t_{\max}^+, ~ t_{\max}^-) \in \mathbb{Z}_+^2$
defined below, and
various maximal power vectors.
\begin{enumerate}
\item\label{step:iterationlengthfuturepastoutputs}
{\em Construct the vectors of the past and of the future time series.}
Take a time $t \in T = \{ 1,2, \ldots, t_1\}$ less or
equal to $t_1/2$.
Denote the length of the tuple of the future and of the past output
time series respectively 
by $(t^+,~ t^-) \in \mathbb{Z}_+$
and set their extrema
such that $(t+t_{\max}^+),(t-t_{\max}^-) \in T$.
\par
Iterate 
from Step \ref{step:setpowermatrices}
to Step \ref{step:lkreduction}
in a Levinson-like manner
by increasing the horizon lengths $(t^+, t^-)$
from the values $(t_{\min}^+, ~ t_{\min}^-)$
to the values
$(t_{\max}^+, ~ t_{\max}^-) \in \mathbb{Z}_+^2$.
\par
Construct the symbolic vectors of 
the future and the past series,
and their values for the each of the time series,
\begin{eqnarray}
  &   & d_{y^+} = t^+ d_y, ~
        d_{y^-} = t^- d_y \in \mathbb{Z}_+, ~ \nonumber \\
      \lefteqn{
        (y^+(t), y^-(t))
        \in (\mathbb{R}^{d_{y^+}} \times \mathbb{R}^{d_{y^-}}) 
      } \label{eq:yplusymin} \\
  & = & \left(
        \left(
        \begin{array}{l}
          y(t+t^+-1) \\
          y(t+t^+-2) \\
          \vdots \\
          y(t) 
        \end{array}
        \right), ~
        \left(
        \begin{array}{l}
          y(t-1) \\
          y(t-2) \\
          \vdots \\
          y(t-t^-)
        \end{array}
        \right)
        \right). \nonumber  \\
      \lefteqn{
        (y^+(t,k), y^-(t,k))
        \in (\mathbb{R}^{d_{y^+}} \times \mathbb{R}^{d_{y^-}}), ~
        k \in \mathbb{Z}_s. \nonumber
      } 
\end{eqnarray}
\item\label{step:setpowermatrices}
{\em Define the power matrices of the future and past
output time series}.
\begin{eqnarray}
  &   & d_{v^+} = d_{y^+}, ~ 
        K_{v^+} = I_{d_{v^+}}
        \in \mathbb{N}^{d_{v^+} \times d_{y^+}}, \nonumber \\
  &   & k_{\max,y} \in \mathbb{N}^{d_{y}}, 
        k_y^* = \max_{i \in \mathbb{Z}_{d_y}} k_{\max,y}(i), \\
      k_{\max,y^-} 
  & = & \left(
        \begin{array}{lll}
          k_{\max,y}^T, & \ldots, & k_{\max,y}^T 
        \end{array}
        \right)^T \in \mathbb{N}^{d_{y^-}}, \nonumber  
\end{eqnarray}
\begin{eqnarray}
      d_{v^-}
  & \leq & ( \prod_{i=1}^{d_y} (k_{\max,y}(i) + 1))^{t^-} 
           \leq (k_y^*+1)^{d_{y^-}}, \nonumber \\
  &   & K_{v^-} \in \mathbb{N}^{d_{v^-} \times d_{y^-}}(k_{\max,y^-}), 
          \nonumber \\
  &   & (L_{v^+},K_{v^+}) = (I_{d_{v^+}}, I_{d_{v^+}}), \nonumber \\
  &   & (L_{v^-}, K_{v^-}) = (I_{d_{v^-}}, K_{v^-}).  \nonumber
\end{eqnarray}
\item\label{step:iterationmonomialvectors}
{\em Iteration of blocks of the power matrix $K_{v^-}$}.
If the row dimension of 
the bounded power matrix $K_{v^-}$ 
of the past outputs
is very high, say larger than 500, 
then execute this step.
Partition the full power matrix $K_{v^-}$
into a finite number of row blocks
$K_{v^-}^{(1)}, \ldots, K_{v^-}^{(m_1)}$
starting with lowest power blocks.
Denote the corresponding row dimensions by
$d_{v^-}^{(1)}, \ldots, d_{v^-}^{(m_1)}$.
\par
Iterate from Step \ref{step:constructionmonomialoutputvectors}
to Step \ref{step:monomialapproximationfuturepast}.
Start with the first block
$(d_{v^-}^{(1)}, L_{v^-}^{(1)}, K_{v^-}^{(1)})$.
After each cycle 
add to the current generator set
indexed by 
$(L_{v^-}^{(g)}, K_{v^-}^{(g)})$
the next block of the power matrix.
\item\label{step:constructionmonomialoutputvectors}
{\em Construct the monomial vectors of the future and the past outputs.}
Construct next for the parameters \\
$(d_{v^+},K_{v^+}, d_{v^-}^{(m)}, K_{v^-}^{(m)})$
set in 
Step \ref{step:setpowermatrices} or 
Step \ref{step:iterationmonomialvectors},
the associated monomial vectors
according to Def. \ref{def:powerboundedmonomialvector},
\begin{eqnarray}
      v^+ (k)
  & = & v(y^+(t,k),d_{y^+},K_{v^+}) 
        = y^+(k)  \in \mathbb{R}^{d_{v^+}},  \\
      V^+(t)
  & = & \left(
        \begin{array}{llll}
          v^+(1) & v^+(2) & \ldots & v^+(s) 
        \end{array}
	\right) \in \mathbb{R}^{d_{v^+} \times s},  \\
      v^- (k)
  & = & v(y^-(t,k),d_{y^-},K_{v^-}^{(m)}) \in \mathbb{R}^{d_{v^-}^{(m)}}, \\
      v_i^-(k)
  & = & \prod_{j=1}^{d_{y^-}} y_j^-(t,k)^{K_{v^-}^{(m)}(i,j)}, ~
        \mbox{see (\ref{eq:vxk})}, \\
      V^-(t)
  & = & \left(
        \begin{array}{llll}
          v^-(1) & v^-(2) & \ldots & v^-(s) 
        \end{array}
        \right) \in \mathbb{R}^{d_{v^-}^{(m)} \times s}.
\end{eqnarray}
Note that the dimension and hence the complexity of $v^-$
is exponential in terms of $d_{v^-} = t^- d_y$. 
\item\label{step:monomialapproximationfuturepast}
{\em Reduce the generator set by (1) linear dependence}.
Compute,
according to Algorithm \ref{alg:approximationpolynomialmapbylinpolyn},
the approximate monomial equation
of the future and the past time series.
\begin{eqnarray}
  &   & (n_1,D_{n_1},C_{v^+},L_1,X,H^*,table_1) \\
  & = & \svdtruncation (dy^+,dy^-,dv_{y^+},dv_{y^-},s, \nonumber \\
  &   & V^+(t),V^-(t),r_1); ~
        H^* \in \mathbb{R}^{d_{v^+} \times d_{v^-}}, \nonumber  \\
      V^+(t)
  & \approx & H^*(t) V^-(t) 
        = C_{v^+} L_1 (y^-(t,*))^{K_{v^-}^{(m)}},
        \label{eq:vplust} \\
  &   & C_{v^+} \in \mathbb{R}^{d_{v^+} \times n_1}, ~
        L_1 \in \mathbb{R}^{n_1 \times d_{v^-}^{(m)}}. 
\end{eqnarray}
\item\label{step:lkreduction}
{\em LK-Reduction}.
Reduce the generator set further
for the matrix tuple $(L_1,K_{v^-}^{(i)})$ 
by deleting those columns of the matrix $L_1$
and the corresponding rows of the matrix $K_1$
whose $l_1$-norm of the column is less than
$r_4 \in (0,1)$ times the $1$-norm of $L_1$.
Thus delete column $j$ of $L_1$ if,
\begin{eqnarray}
  &   & \sum_{i=1}^n |L_1(i,j)| 
        \leq r * \|L_1\|_{L_1}
        = r * \max_{j \in \mathbb{Z}_{n}}
                \sum_{i=1}^n |L_1(i,j)|, \nonumber \\
  &   & (L_{v^-}^{(g)},K_{v^-}^{(g)}), ~
        \mbox{result.} \label{eq:lkreduced}
\end{eqnarray}
\item\label{step:transcedencebasis}
{\em Reduce the generator set by 
(2) elimination of products of generators}.
Compute a new generator set 
with possibly less generators
according to Step 2 of the algorithm
of Def. \ref{algorithm:approximatepolynomialfactorization}.
Starting from equation
(\ref{eq:lkreduced}).
with $(C_{v^+}, (L_{v^-}^{(g)},K_{v^-}^{(g)}))$. 
The result is, 
\begin{eqnarray}
      y^+(t)
  & = & v^+(t) 
        \approx C_{v^+} L_{v^-}^{(g)} (y^-(t))^{K_{v^-}^{(g)}}  \nonumber \\
  & \approx & h_{io}^+(x(t)) 
        = L_{h_{io}^+} x(t)^{K_{h_{io}^+}}, ~
        \mbox{(\ref{eq:vxk})} \\
      x(t)
  & = & g_{io}(v^-(t)) = L_{g_{io}} y^-(t)^{K_{g_{io}}}
        \in \mathbb{R}^{n}.
\end{eqnarray}
\item\label{step:computeoutputequation}
{\em Compute the output equation.}
\begin{eqnarray*}
      y(t)
  & = & P_{y(t)} y^+(t)
        \approx
        P_{y(t)} h_{io}^+(x(t)) 
        = h_o(x(t)), \label{eq:outputeqcomputed} \\
  &   & P_{y(t)} \in \mathbb{R}^{d_y \times d_{y^+}}, ~
        \mbox{a projection.} \nonumber
\end{eqnarray*}
The next steps aim at the computation of the
system dynamics, see equation
(\ref{eq:systemdynamicspolynomial}).
\item
{\em Compute the value of the next state.}
First compute 
the past time series at the next time index $(t+1)$.
Secondly, compute the value of $X(t+1)$
for each time series.
\begin{small}
\begin{eqnarray}
      y^-(t+1,j)
  & = & \left(
        \begin{array}{l}
          y(t,j) \\ y(t-1,j) \\ \vdots \\ y(t-t^- +1,j) 
        \end{array}
        \right) \in \mathbb{R}^{d_{y^-}}, \\
      v^-(t+1,j)
  & = & v(y^-(t+1,j),d_{y^-},K_{g_{io}}) \in \mathbb{R}^{d_{v^-}}, \\
      V^-(t+1)
  & = & \left(
        \begin{array}{llll}
          v^-(t+1,1) & v(t+1,2) & \ldots & v(t+1,s)
        \end{array}
        \right) \\
  &   & \in \mathbb{R}^{d_{v^-} \times s}; ~ 
        d_{v_{x}} = n \in \mathbb{Z}_+, \nonumber \\
      V_{x}(t+1)
  & = & X(t+1)
        = L_{g_{io}} V^-(t+1) \in \mathbb{R}^{n \times s}. \nonumber
\end{eqnarray}
\end{small}
\item\label{step:computemonomialxy}
{\em Compute the monomial vector
of the current state and the current output.}
\begin{small}
\begin{eqnarray}
      d_{(x,y)} 
  & = & n + d_y, ~ d_{v_{(x,y)}} \in \mathbb{Z}_+, \\
      d_{v_{(x,y)}}
  & = & \prod_{i=1}^n [k_{\max,x}(i)+1]
        \prod_{j=1}^{d_y} [k_{\max y2}(j) + 1], \\
  &   & \mbox{choose} ~ 
        k_{\max,x} \in \mathbb{N}^n, ~
        k_{\max,y2} \in \mathbb{N}^{dy}, \nonumber \\
  &   & k_{\max,(x,y)} =
        (
        \begin{array}{ll}
          k_{\max,x}^T 
          k_{\max,y2} ^T
        \end{array}
        )^T \in \mathbb{N}^{d_{(x,y)}}, \nonumber \\
  &   & K_{v_{(x,y)}} \in 
          \mathbb{N}^{d_{v_{(x,y)}} \times d_{(x,y)}}
          (k_{\max,(x,y)}), \nonumber \\      
      v_{(x(t),y(t))}(k)
  & = & v(
        \left(
        \begin{array}{l}
          x(t,k) \\
          y(t,k)
        \end{array}
        \right),
        (n+d_y), K_{v_{(x,y)}}
        ) \\
      V_{x,y} 
  & = & \left(
        \begin{array}{lll}
          v_{(x(t),y(t))}(1) & \ldots & v_{(x(t),y(t))}(s) 
        \end{array}
        \right). 
\end{eqnarray}
\end{small}
Note that the complexity of the expression of
$d_{v_{x,y}}$ 
is exponential in $n+d_y$. 
\item
{\em Reduce the generator set by (1) linear dependence}.
Compute the approximate polynomial function of the next future state
depending on the vector of the current state and of the current output.
Compute according to 
Def. \ref{alg:approximationpolynomialmapbylinpolyn},
\begin{eqnarray}
  &   & (n_2,D_{n_2},C_2,L_2,X_2,H_2^*,table_2) \\
  & = & \svdtruncation 
        (n,d_{(x,y)},n,dv_{(x,y)},s,V_x(t+1), \nonumber \\
  &   & V_{x,y},r_2); ~
        L_3 = H_2^*  \in \mathbb{R}^{n \times d_{v_{(x,y)}}}. 
        \nonumber \\
      X(t+1) 
  & = & V_{x(t+1)} 
        \approx L_3 V_{(x(t),y(t))}, \\
      x(t+1) 
  & = & L_3 ~ (x(t),y(t))^{K_3}.
\end{eqnarray}
\item\label{step:approximatemonomialmap}
{\em Approximate the monomial map.}
Reduce the matrix tuple $(L_3,~ K_3)$
to $(L_{f_{io}},K_{f_{io}})$
as in Step \ref{step:lkreduction}
for the matrix pair $(L_2, K_{v^-}^{(i)})$.
\item
{\em Compute the polynomial system.}
Finally compute the discrete-time polynomial system without input 
in the observer form,
by writing the linear map of nomomials
as a vectorial polynomial function,
\begin{eqnarray}
      f_{o}(x,y)
  & = & L_{f_{io}} (x,y)^{K_{f_{io}}}, ~
        \mbox{see (\ref{eq:vxk}),} \\
      x_o(t+1)
  & = & f_o(x_o(t),y(t)), ~ x_o(0) = x_{o,0}, \\
      y_o(t|t-1)
  & = & h_o(x_o(t)), ~
        \mbox{from Step (\ref{step:computeoutputequation})},\\
     X_{o,0}
  & = & X(t) \in \mathbb{R}^{n \times s}.
\end{eqnarray}
Output $(n, ~ f_o, ~ h_o, ~ \{ X_{o,0}(k), ~ y^+(k), ~ y^-(k), ~ k \in \mathbb{Z}_s\})$.
\end{enumerate}
\end{definition}
The user is advised 
to take $t^+$ and $t^-$ equal to about 4 times 
the dimension of the expected state space. 
\section{Examples}\label{sec:examples}
\begin{example}\label{ex:example1}
For a set of times series generated
by a simple polynomial system
with a one-dimensional output and 
a two-dimensional state vector,
a computer program for 
Algorithm~\ref{alg:siapolsysdtnoinput}
has computed the following observer polynomial system.
\begin{eqnarray*}
      x(t+1)
  & = & f_o(x(t),y(t))
        = L_{f_o} (x_o(t),y(t))^{K_{f_o}}, \\
      y(t|t-1)
  & = & h_o(x(t)) = C x_o(t), ~
      C
      = \left(
        \begin{array}{ll}
          -0.0225 & 0.0336
        \end{array}
        \right), \\
      L_{f_o}
  & = & \left(
        \begin{array}{rrrrrr}
          0.009 & .089 &  0.023 & .571 & -.004 & -0.020 \\
         -0.015 & .309 & -0.014 & .074 &  .008 &  0.212
        \end{array}
        \right) \nonumber \\
      (x_o,y)^{K_{f_{io}}}
  & = & (
         x_{o,1} x_{o,2} y,
         x_{o,1} x_{o,2},
         x_{o,1} y,
         x_{o,1},
         x_{o,2} y,
         x_{o,2}
        )^T. 
\end{eqnarray*}
\end{example}
\begin{small}

\end{small}
\appendix
\section{Notation}\label{ap:notation}
The set of the integers is denoted by
$\mathbb{Z}$ and the positive integers by
$\mathbb{Z}_+ = \{ 1, 2, \ldots \}$.
The set of the natural numbers is denoted by
$\mathbb{N} = \{ 0, 1, 2, \ldots \}$
and by $\mathbb{N}^n$ its $n$-tuples.
For any $n \in \mathbb{Z}_+$,
denote $\mathbb{Z}_n = \{ 1, 2, \ldots, n\}$.
The set of the real numbers is denoted by $\mathbb{R}$,
the positive real numbers by $\mathbb{R}_+ = [0,\infty)$, and
the strictly-positive real numbers by $\mathbb{R}_{s+} = (0,\infty)$.
The vector space of $n$-tuples of the real numbers
is denoted by $\mathbb{R}^n$, for $n \in \mathbb{Z}_+$.
The set of matrices with entries in the real numbers 
of size $n \times m$, for $n, m \in \mathbb{Z}_+$,
is denoted by $\mathbb{R}^{n \times m}$.
A diagonal matrix of the set of square real matrices
$\mathbb{R}^{n \times n}$ for $n \in \mathbb{Z}_+$
is a matrix such that its off-diagonal elements all zero
and the set of such matrices is denoted by $\dnn$.
The subset of positive diagonal matrices 
is defined by the condition 
that $D_{i,i} \geq 0$ for all $i \in \mathbb{Z}_n$ and 
it is denoted by $\dposnn$.
Similarly, $\dsposnn$.
\begin{definition}\label{def:diagmattruncation}
The {\em truncation operation} of a positive diagonal matrix based
on the $l_1$-norm of the diagonal. \\
{\em Data.}
$(n, D, r) \in (\mathbb{Z}_+ \times \dposnn \times (0,1))$,
$D \neq 0$,
where $r$ is an approximation threshold. 
Assume that 
$D_{1,1} \geq D_{2,2} \geq \ldots \geq D_{n,n} \geq 0$.
\begin{enumerate}
\item
Compute the $l_1$-norm of the diagonal elements
of the diagonal matrix $D$,
$\|\diag (D)\|_{l_1} = \sum_{i=1}^n D_{i,i}$.
\item
Compute\\
$n_r = \argmin_{j \in \mathbb{Z}_n}
        \{ \sum_{i=1}^{j} D_{i,i} 
           / \|\diag (D)\|_{l_1} 
           \geq r 
        \}. 
$
\item
Construct the approximant positive diagonal matrix
\begin{eqnarray}
      D_{n_r,i,i} 
  & = & D_{i,i}, \forall ~ i \in \{ 1,2, \ldots, n_r \}, ~
        D_{n_r} \in \mathbb{R}_{{\rm diag,+}}^{n_r \times n_r}, 
          \nonumber \\
      D_{r}
  & = & \left( 
        \begin{array}{ll}
          D_{n_r} & 0 \\
          0   & 0
        \end{array}
        \right) \in \dposnn.
        \label{eq:dr}
\end{eqnarray}
\item
Output 
$(n_r, D_r, \Table_1) \in 
(\mathbb{N} \times \dposnn \times \mathbb{R}^{n \times 2})$,
where \\
$\Table_1 = \{ (j, \sum_{i=1}^j D_{i,i}/ \|\diag (D)\|_{l_1}), ~ j \in \mathbb{Z}_n \}$.
\end{enumerate}
\end{definition}
\section{Monomials and Monomial Vectors}\label{ap:monomialvectors}
\begin{definition}\label{def:monomialvectors}
Consider a set of $n \in \mathbb{Z}_+$ commutative variables 
denoted by $x = (x_1, \ldots, x_n)$.
A {\em monomial} 
is a term of a polynomial 
defined by the formulas,
\begin{eqnarray*}
      x^k 
  & = & \prod_{i=1}^n x_i^{k(i)}
        = x_1^{k(1)} x_2^{k(2)} \ldots x_n^{k(n)}, ~
        k \in \mathbb{N}^n, \\
      \gmon [x]
  & = & \{ x^k ~ | ~ \forall ~ k \in \mathbb{N}^n \}, ~
        A(\gmon)[x] = A_{\mathbb{R}^n}[x], \nonumber \\ 
      p(x) 
  & = & \sum_{k \in \mathbb{N}^n} c(k) x^k \in \mathbb{R}[x], ~
        \forall ~ k \in \mathbb{N}^n, ~ c(k) \in \mathbb{R}.
\end{eqnarray*}
Call 
$x^k$ a {\em monomial} in the indeterminates $x$,
a vector $k \in \mathbb{N}^n$ a {\em power vector},
$\mathbb{N}^n$ a {\em set of power vectors},
$\gmon[x]$ the set of all monomials of $x$, and
$p \in \mathbb{R}[x]$ a 
{\em polynomial in monomial representation}.
\end{definition}
\par
There exists a bijective correspondence between 
the set of power vectors $\mathbb{N}^n$ 
and the set $\gmon [x]$ of monomials 
in the indeterminates $x$,
described by the map $k \mapsto x^k$.
The partially-ordered set $\mathbb{N}^n$ 
may be equiped with a monomial ordering.
Below the specific monomial ordering
called the {\em lexicographic order relation} on
the index set of power vectors $\mathbb{N}^n$ 
will be used, see \cite[Def. 2.2.3]{cox:little:oshea:1992}.
It is denoted by $\orderlexg$.
By the bijective correspondence between
$\mathbb{N}^n$ and $\gmon [x]$
the lexicographic ordering of $\mathbb{N}^n$
is transformed into a 
lexicographic ordering on $\gmon [x]$
which is again denoted by $\orderlexg$ 
and which will be called 
the {\em lexicographic ordering} of $\gmon [x]$.
\begin{definition}\label{def:powerboundedmonomialvector}
Define the {\em power-bounded monomial vector}
of a set of $n \in \mathbb{Z}_+$ commutative variables 
by the following formulas.
\begin{eqnarray}
  &   & k_{\max} \in \mathbb{N}^n, \nonumber \\
      \mathbb{N}^n(k_{\max})
  & = & \{ k \in \mathbb{N}^n |
           0 \leq k(i) \leq k_{\max}(i), ~ \forall ~ i \in \mathbb{Z}_n
        \}, \nonumber \\
      k^*
  & = & |\mathbb{N}^n(k_{\max})| 
        =  \prod_{i=1}^n (k_{\max}(i) + 1) \in \mathbb{Z}_+,  \nonumber \\
  &   & \mbox{choose} ~ 
        d_v \in \mathbb{Z}_+, ~
        0 < d_v \leq k^*, ~ \nonumber \\
      \mathbb{N}^{d_v \times n}(k_{\max})
  & = & \{ K \in \mathbb{N}^{d_v \times n} |
           0 \leq K(i,j) \leq k_{\max}(j)
        \},  \nonumber \\
  &   & \mbox{choose} ~ 
        K_v \in \mathbb{N}^{d_v \times n} (k_{\max}), 
        \mbox{and define,}  \nonumber \\
      v 
  & = & v(x,n,K_v) = x^{K_v} \in \mathbb{R}^{d_v}, \label{eq:vxk} \\
      v(x,n,K_v)_i
  & = & x_1^{K_v(i,1)} 
          x_2^{K_v(i,2)} 
          \ldots
          x_n^{K_v(i,n)}; ~
        h \in \mathbb{R}^{d_v}, \nonumber \\
      p(x) 
  & = & \sum_{k \in K_v} c(k) x^k 
        = h^T v(x,n,K_v)
        = h^T x^{K_v}. \nonumber
\end{eqnarray}
Call
$k_{\max}$ the {\em maximal power vector} and
$k_{\max}(i)$ the {\em maximal power} of $x_i$; 
$\mathbb{N}^n(k_{\max})$
the {\em bounded power vector set},
which set inherits the lexicographic ordering
of the elements of $\mathbb{N}^n$;
$K_v \in \mathbb{N}^{d_v \times n}(k_{\max})$
the {\em power matrix} of the monomial vector $v(x,n,K_v)$,
where the row elements of $K_v$ are
the powers vectors of the vector $v$ 
in decreasing lexicographic order;
$v(x,n,K_v)$ a (symbolic) {\em monomial vector},
which contains all monomials indexed by 
$K_v \in \mathbb{N}^{d_v \times n}(k_{\max})$
in their lexicographic order from high to low order;
the number of components of this monomial vector equals
$d_v \leq k^*$; and
finally call the equation
$p(x) =  h^T v(x,n,K_v) = h^T x^{K_v}$, 
the {\em power-bounded monomial representation}
of the polynomial $p$.
\end{definition}
\begin{example}
There follows an example of a monomial vector. 
\begin{eqnarray*}
      x
  & = & (x_1, ~ x_2), ~ n = 2, ~ k_{\max} = (2,1)^T, ~ \\
  &   & d_v = k^* = (2+1) \times (1+1) = 6, \\
      v(x,2,(2,1))
  & = & \left(
        \begin{array}{l}
          x_1^2 x_2 \\
          x_1^2 \\
          x_1 x_2 \\
          x_1 \\
          x_2 \\
          1
        \end{array}
        \right); ~~
      K_v
      = \left(
        \begin{array}{ll}
          2 & 1 \\   
          2 & 0 \\
          1 & 1 \\
          1 & 0 \\
          0 & 1 \\
          0 & 0
        \end{array}
        \right) \in \mathbb{N}^{6 \times 2}(
        \left(
        \begin{array}{l}
          2 \\ 1
        \end{array}
        \right) ).
\end{eqnarray*}
\end{example}
\section{Polynomials}\label{ap:polynomials}
In this appendix and the next one
several aspects of the commutative algebra
of polynomials maps are described
including algebraic geometry.
See
\cite{becker:weispfenning:1993,demoor:2014},
\cite{cox:little:oshea:1992}, and
\cite{zariski:samuel:1958}.
\par
Let $n \in \mathbb{Z}_+$ be a positive integer.
The {\em ring of polynomials} in $n$ variables
with real coefficients is denoted
by $\mathbb{R}[x_1,\ldots,x_n]$.
The simplified notation of
$\mathbb{R}[x]$ will be used if it is understood
that the variable $x$ has $n$ components.
Examples are
$2 x^2 + 3 x + 4 \in \mathbb{R}[x]$ and
$21 x_1^2 x_2 + 11 x_1 x_2 + 1 x_2 \in \mathbb{R}[x_1,x_2]$.
\par
Below polynomial functions 
of tuples of the real numbers $X = \mathbb{R}^n$
are needed.
A {\em polynomial function} on $\mathbb{R}^n$
is a map $p: X \rightarrow \mathbb{R}^n$ 
for which there exists a set of polynomials 
$q_1, \ldots, q_n \in \mathbb{R}[x_1, x_2, \ldots, x_n]$ such that
$p_i = q_i$ on $\mathbb{R}^n$ for all $i \in \zn$.
Denote by $A_X$ the set of all polynomials 
on $X = \mathbb{R}^n$.
\begin{definition}
The {\em monomial representation}
of a finite set of polynomials
in the indeterminates $x = x_1 \ldots x_{d_x}$
power bounded by the vector $k_{\max} \in \mathbb{N}^n$,
where the polynomials are the components of $L x^K$,
is defined by the notation,
\begin{eqnarray}
      G 
  & = & \left\{ 
        \begin{array}{l}
           L x^K \in \mathbb{R}[x]| 
           (L,K) \\
           \in 
           ( \mathbb{R}^{\card_G \times d_v}
             \times
             \mathbb{N}^{d_v \times d_x}(k_{\max})
           )
        \end{array}
        \right\}.   \label{eq:generatorset}
\end{eqnarray}
\end{definition}
\begin{example}
Consider the set of polynomials, 
\begin{eqnarray*}
      L 
  & = & \left(
        \begin{array}{ll}
          0.1 & 0.2 \\
          0.3 & 0.4
        \end{array}
        \right), ~
      K
      = \left(
        \begin{array}{ll}
           3 & 1 \\
           1 & 2
        \end{array}
        \right), ~
        x^K = 
        \left(
        \begin{array}{l}
          x_1^3 x_2 \\
          x_1 x_2^2
        \end{array}
        \right), \\ 
  &   & \left\{ 
        p(x) =
        L x^K =
        \left(
        \begin{array}{l}
          0.1 x_1^3 x_2 + 0.2 x_1 x_2^2 \\
          0.3 x_1^3 x_2 + 0.4 x_1 x_2^2
        \end{array}
        \right)
        \right\}.
\end{eqnarray*}
\end{example}
\section{A Set of Generators}
A {\em subalgebra} $A_1$ of the algebra $A_X$
is a subset $A_1 \subseteq A_X$
such that the algebraic operations
of $A_1$ are obtained from those of $A_X$ by restriction
and such that it is also an algebra itself 
in terms of those operations.
\par
Consider a subset $G \subseteq A_X$.
The smallest subalgebra of $A_X$
containing $G$ exists, 
it is called the {\em algebra generated} by the set $G$,
and it is denoted by $A_X(G) \subseteq A_X$.
A subalgebra $A_1 \subseteq A_X$ is called
{\em finitely generated} 
if there exists a finite subset $G_f \subset A_X$ 
such that $A_1 = A_X (G_f)$.
\begin{definition}
Call the finite set $G$,
see (\ref{eq:generatorset}),
a {\em generator set} 
of the subalgebra $A \subseteq \mathbb{R}[x]$
and any row component of $L x^K \in G$ a {\em generator} of $A$
if $A = A(G)$.
Call it a {\em nontrivial generator set}
if no column of the matrix $L$ in the representation
is entirely zero.
Call $G$ a {\em minimal generator set}
of the algebra $A = A(G)$
if it is nontrival 
and for any other generator set $H$
it holds that $\card_G \leq \card_H$.
A {\em set of generators} 
of a {\em finite set} of polynomials $H$
is a finite set of polynomials $G$ such that
\begin{eqnarray}
      A(H)
  & = & A(G), ~
      \mbox{where,} ~
      H = \{ p_1, \ldots, p_{\card_H} \in \mathbb{R}[x] \}, \nonumber \\
      G 
  & = & \{ g_1, \ldots, g_{\card_G} \in \mathbb{R}[x] \}. \nonumber
\end{eqnarray}
\end{definition}
\par
The set of real numbers is also a ring.
A finite subset $\{p_1, \ldots, p_k\} \subset \mathbb{R}[x]$ is 
called {\em algebraically dependent} over $\mathbb{R}$
if there exists a nonzero polynomial $f \in \mathbb{R}[p]$
such that $f(p_1, \ldots, p_k) = 0$.
It is called {\em transcedental} otherwise.
See \cite[I, \S 17, p.28]{zariski:samuel:1958}.
An extension not described here because of lack of space
is to define a transcedence basis for the algebraic structure used
which then allows the use of the algorithms of
\cite{muellerquade:steinwandt:2000} for the computation
of such a basis.
\begin{problem}
Consider a finite set of polynomials
$H \subset \mathbb{R}[x]$.
Construct a minimal set of generators
$G \subset \mathbb{R}[x]$ of $H$. 
\end{problem}
The above problem is {\em equivalent} 
to the problem of constructing 
a polynomial factorization 
of a polynomial map as defined next.
\begin{definition}
A {\em polynomial factorization}
of a polynomial map $y = f(u)$ with 
$d_y, ~ d_u \in \mathbb{Z}_+$,
$f: \mathbb{R}^{d_u} \rightarrow \mathbb{R}^{d_y}$,
is defined to be a factorization of the form,
\begin{eqnarray}
      y
  & = & f(u) = h(g(u)) = h(x), ~
        h \in \mathbb{R}[x], \\
      x 
  & = & g(u) \in X \subseteq \mathbb{R}^{d_x}, ~
        g \in \mathbb{R}[u], \\
      G_f
  & = & \{ f_1, \ldots, f_{d_y} \in \mathbb{R}[u] \}, ~
        G_g = \{ g_1, \ldots, g_{d_x} \in \mathbb{R}[u] \}, 
        \nonumber \\
      A(G_f)
  & = & A(G_g),
\end{eqnarray}
hence $G_g$ is a set of generators of $A(G_f)$.
Note that $G_g$ is a minimal set of generators
if $d_x \in \mathbb{Z}_+$ is minimal
over all factorizations.
\end{definition}
\section{Approximation of a Generator Set}\label{ap:approximationgenset}
\begin{problem}\label{problem:approximatepolynomialfactorization}
The problem of {\em approximate polynomial factorization}.
Consider a polynomial function
$y = f(u)$ as defined in Section \ref{ap:polynomials}.
(The notation of this section differs from the main body of the paper.)
Determine an approximate polynomial factorization
of the form,
\begin{eqnarray}
     y
  & = & f(u) \approx h(g(u)) = h(x), ~ x = g(u), \\
  &   & d_x \in \mathbb{Z}_+, ~
        g \in \mathbb{R}[u], ~
        h \in \mathbb{R}[x]. \nonumber
\end{eqnarray}
The approximation criterion of the expression
$[ f(u) - h(g(u))]$ is not specified.
\end{problem}
\begin{definition}\label{algorithm:approximatepolynomialfactorization}
{\em The approximate polynomial factorization algorithm}.
Data.
$y = f(u)$, $d_y, ~ d_u \in \mathbb{Z}_+$,
$f \in \mathbb{R}[u]$.
\begin{enumerate}
\item
Construct the approximation consisting of a 
linear\---po\-ly\-no\-mi\-al factorization by
the algorithm of Def. \ref{alg:approximationpolynomialmapbylinpolyn},
\begin{eqnarray}
     y
  & = & f(u) \approx C_r g_r(u) = C_r x_r, ~ x_r = g_r(u), \\
  &   & d_{x_r} \in \mathbb{Z}_+, ~
        x_r \in \mathbb{R}^{d_{x_r}}, ~
        C_r \in \mathbb{R}^{d_y \times d_{x_r}}.
\end{eqnarray}
\item
For the linear--polynomial factorization of Step 1 
construct a new generator set of, possibly, 
lower cardinality than before,
or, equivalently, a polynomial factorization,
\begin{eqnarray}
      y
  & = & f(u) \approx C_r g_r(u) = h(g(u)) = h(x), \\
      x
  & = & g(u),  \\
  &   & d_x \in \mathbb{Z}_+, ~
        x \in \mathbb{R}^{d_x}, ~
        h \in \mathbb{R}[x], ~
        g \in \mathbb{R}[u], \nonumber
\end{eqnarray}
by polynomial factorizations of the components
of $g_r$.
For example, if for polynomial $g_{r,m}$
there exists a factorization of the form
$g_{r,m} = g_{r,i} g_{r,j} + g_{r,s}$, 
where $g_{r,i}, ~ g_{r,j}, ~g_{r,s}$ have lower
powers than those of $g_{r,m}$,
see \cite[Sec. 5.1.]{becker:weisspfenning:1993}.
\end{enumerate}
\end{definition}
Comments follow.
(1) A Gr\"{o}bner basis algorithm is not appropriate
for Problem \ref{problem:approximatepolynomialfactorization}
of the polynomial function $y = f(u)$
because that would not be an approximation.
In addition, the computational complexity is too high.
(2) A Gr\"{o}bner basis algorithm may well 
be useful for Step 2 of 
Algorithm \ref{algorithm:approximatepolynomialfactorization}.
For the computations a simple procedure is used for Step 2 of 
Algorithm \ref{algorithm:approximatepolynomialfactorization}.
The procedure is not described here in detail
because of lack of space.
\par
The transformation for a polynomial factorization
is briefly described.
Note that if
\begin{eqnarray}
     y 
  & = & f(u) \approx C_r g_r(u), ~ x = g_r(u), \\
      g_{r,m}(u)
  & = & x_{r,m} = x_{r,i} x_{r,j} = g_{r,i}(u) g_{r,j}(u), ~
        \mbox{then,} \\
      x_r
  & = & \left(
        \begin{array}{llllll}
          x_1 & \ldots & x_{m-1} & x_{m+1} & \ldots & x_{d_x} 
        \end{array}
        \right)^T,  \nonumber \\
      x 
  & = & \left(
        \begin{array}{lllllll}
          x_1 & \ldots & x_{m-1} & x_i x_j & x_{m+1} & \ldots & x_{d_x} 
        \end{array}
        \right)^T, \nonumber \\
  & = & P_x
        \left(
        \begin{array}{llll}
          x_i x_j & x_1 & \ldots & x_{d_x} 
        \end{array}
        \right) ^T  
        = P_x v(x_r) = P_x x_r^{K_{x_r}}, \nonumber \\
      y
  & \approx & C_r x = C_r P_x v(x_r) = C v(x_r) = h(x_r), \\
      g(u)
  & = & P_{g_r} g_r(u), \\
      y
  & \approx & C_r x = C_r g_r(u) = h(x_r), ~
        x_r = g(u),
\end{eqnarray}
It is not claimed that the above procedure
determines a minimal generator set which is not
true in general.
\begin{definition}\label{def:monomialequationdatamatrices}
{\em The monomial equation of data matrices}.
Consider the polynomial equation $y = f(u)$ and
its monomial representation.
Consider the case in which one is provided
several tuples of values of input and output vectors,
$\{(\overline{y}_i,\overline{u}_i) \in Y \times U | ~ i = 1, 2, \ldots, s \}$.
A {\em monomial equation of the data matrices}
for this set of tuples 
is then a linear map 
represented by the coefficient matrix $H$.
\begin{eqnarray}
      V_y
  & = & H V_u,  ~
        H \in \mathbb{R}^{d_{v_y} \times d_{v_u}}, \\
      V_y
  & = & \left(
        \begin{array}{llll}
            v(\overline{y}_1,d_y,K_{v_y})  
          & \ldots 
          & v(\overline{y}_s,d_y,K_{v_y}) 
        \end{array}
        \right) 
        \in \mathbb{R}^{d_{v_y} \times s}, \nonumber \\
      V_u
  & = & \left(
        \begin{array}{llll}
            v(\overline{u}_1,d_u,K_{v_u})  
          & \ldots 
	  & v(\overline{u}_s,d_u,K_{v_u}) 
	  \end{array}
        \right) \in \mathbb{R}^{d_{v_u} \times s}. \nonumber 
\end{eqnarray}
\end{definition}
\begin{definition}\label{alg:approximationpolynomialmapbylinpolyn}
{\em Linear approximation of a polynomial map}.\\
\cite[Sec. 6.1]{golub:vanloan:1983}.\\
This algorithm is called $SVDtrunction$ 
in Steps 5 and 11 
of Algorithm \ref{alg:siapolsysdtnoinput}.
{\em Data.}
$(d_y, d_u, d_{v_y}, d_{v_u}, s, V_y, V_u, r)
\in (\mathbb{Z}_+^5 \times 
\mathbb{R}^{d_{v_y} \times s} \times
\mathbb{R}^{d_{v_u} \times s} \times (0,1))$. 
\begin{enumerate}
\item
Compute the singular value decomposition
of the data matrix of the inputs,
\begin{eqnarray}
     V_u
  & = & V_1^T S V_2 \in \mathbb{R}^{d_{v_u} \times s}, \\
  &   & V_1 \in \mathbb{R}^{d_{v_u} \times d_{v_u}}, ~
        V_2 \in \mathbb{R}^{s \times s}, ~
        \mbox{orthogonal matrices,} \nonumber \\
      S
  & = & \left(
        \begin{array}{ll}
          D & 0 \\
          0 & 0 
        \end{array}
        \right) \in \mathbb{R}^{d_{v_u} \times s}, ~
        n_1 \in \mathbb{N}, \nonumber \\
      D
  & = & \diag (d_1, ~ d_2, \ldots, ~ d_{n_1})
        \in \mathbb{R}_{{\rm \diag,s+}}^{n_1 \times n_1}, 
        \nonumber \\
  &   & d_1 \geq d_2 \geq \ldots \geq d_{n_1} > 0.  \nonumber
\end{eqnarray}
\item
Compute 
according to the algorithm of Def. \ref{def:diagmattruncation},
of the diagonal matrix $D$
its truncation $D_n$ upto the approximation fraction,
\begin{eqnarray}
  &   & (n, D_n, \Table_1) = \matrixdiagtrunc (n_1, D), \\
      S_n^+
  & = & \left(
        \begin{array}{ll}
          D_n^{-1} & 0 \\
          0   & 0
        \end{array}
        \right) \in \mathbb{R}^{s \times d_{v_u}}, ~
        D_n \in \dsposnn, \nonumber \\
      V_u^+
  & = & V_2^T S_n^+ V_1 \in \mathbb{R}^{s \times d_{v_u}}, ~
      H^*
      = V_y V_u^+ \in \mathbb{R}^{d_{v_y} \times d_{v_u}}. \nonumber
          \label{eq:wstarmonomialregressiocoeffmatrix}
\end{eqnarray}
\item
Compute the factorization according to,
\begin{eqnarray}
     H^* 
  & = & [ V_y V_2^T
          \left(
          \begin{array}{l}
            D_n^{-1} \\ 0 
          \end{array}
          \right)
        ] 
        [
           \left(
           \begin{array}{ll}
             I_n & 0 
           \end{array}
           \right) V_1
        ] = C L, ~ \\
      L
  & = & \left(
        \begin{array}{ll}
          I_n & 0 
        \end{array}
        \right) V_1 \in \mathbb{R}^{n \times d_{v_u}}, ~
        X = L V_u \in \mathbb{R}^{n \times s}, \nonumber \\
      C
  & = & V_y V_2^T
          \left(
          \begin{array}{l}
            D_n^{-1} \\ 0 
          \end{array}
          \right)
          \in \mathbb{R}^{d_{v_y} \times n}, ~\\
      V_y
  & \approx & H^* V_u = C L V_u = C X.
\end{eqnarray}
\item
Output $(n, ~ D_n, ~ C, ~ L, ~ X, ~ H^*,\Table)$ with \\
$\Table = \{(j, \sum_{i=1}^j D_{i,i}/ \|\diag(D)\|_{l_1} ), ~ j \in \mathbb{Z}_{n_1} \}$.
\end{enumerate}
\end{definition}

\end{document}